\documentclass[11pt,a4paper]{article}
\usepackage{epsfig}
\usepackage{amssymb}
\usepackage{amsthm}
\usepackage{latexsym}
\usepackage{epic}
\usepackage{eepic}
\usepackage{amsmath}
\usepackage{epstopdf}
\usepackage{multirow}
\usepackage[applemac]{inputenc}

\newcommand{\mc}{\mathcal}

\newcommand{\ph}{{\tilde{p}}}
\newcommand{\ep}{\epsilon}
\newcommand{\bep}{b_{\epsilon}(d)}
\newcommand{\Ex}{\mathrm{E}}
\newcommand{\Var}{\mathrm{Var}}
\newcommand{\Cov}{\mathrm{Cov}}
\newcommand{\qt}{\tilde{q}}
\newcommand{\Qt}{\tilde{Q}}
\newcommand{\weak}{\stackrel{\rm D}{\to}}
\newcommand{\ms}{\stackrel{\rm m.s.}{\to}}
\newcommand{\ninf}{\mbox{ as } n \to \infty}
\newcommand{\Or}{{\rm O}}

\newcommand{\oor}{{\rm o}}

\theoremstyle{plain}
\newtheorem{theorem}{Theorem}
\newtheorem{lemma}{Lemma}

\theoremstyle{plain}

\newenvironment{reflist}{\begin{list}{}
         {\itemsep=20pt \parsep=3pt
          \topsep=0pt  \parskip=20pt  \listparindent=-.15in
           \leftmargin= 0.15in }
         \item \ \vspace{-.35in} }
         {\end{list}}

\begin{document}

\author{David K\"{a}llberg,\quad Oleg\ Seleznjev \vspace{0.5cm}\hfill \\
Department of Mathematics and Mathematical Statistics\\
Ume{\aa } University, SE-901 87 Ume\aa , Sweden\\
}
\date{}
\title{Estimation of quadratic density functionals \\ under $m$-dependence}
\maketitle

\begin{abstract}
In this paper, we study estimation of certain integral functionals of one or two densities with samples from stationary $m$-dependent sequences. We consider two types of $U$-statistic estimators for these functionals that are functions of the number of $\ep$-close vector observations in the samples.
We show that the estimators are consistent and obtain their rates of convergence under weak distributional assumptions. In particular, we propose estimators based on incomplete $U$-statistics which have favorable consistency properties even when $m$-dependence is the only dependence condition that can be imposed on the stationary sequences.
The results can be used for divergence and entropy estimation, and thus find many applications in statistics and applied sciences.
\end{abstract}
\medskip

 \baselineskip=3.4 ex



\section{Introduction}
Let $\{X_i\}_{i=1}^\infty$ and $\{Y_i\}_{i=1}^\infty$ be strictly stationary  $m$-dependent sequences of random $d$-vectors, and denote by $p_X(x)$ and $p_Y(x), x \in R^d$, the densities of  $X_1$ and $Y_1$, respectively.
We consider the problem of estimating the (quadratic) integral functional
$$
q_{k,l} = q_{k,l}( p_X,  p_Y)  := \int_{R^d}p_X(x)^k p_Y(x)^l dx, \quad   k = 2-l = 0,1,2,
$$
using samples $\{X_1,\ldots, X_n\}$ and $\{Y_1, \ldots, Y_n\}$ from these sequences.

Several characteristics in mathematical statistics and information theory are expressed in terms of the functionals $\{q_{k,l} \}_{k=0,1,2}$. Examples include measures of divergence between the distributions of $X_1$ and $Y_1$, e.g., 
the integrated squared difference
\begin{equation*}\label{D}
D(p_X, p_Y) := \int_{R^d} (p_X(x)-p_Y(x))^2 \, dx = q_{2,0} - 2q_{1,1} +  q_{0,2},
\end{equation*}
and functions for quantifying uncertainty in $X_1$, e.g.,
the R\'{e}nyi entropy (R\'{e}nyi, 1970)
\begin{equation*}\label{h2}
h_2( p_X) := - \log \left( \int_{R^d} p_X(x)^2 \, dx \right) = -\log q_{2,0}. 
\end{equation*}
The functional $q_{2,0}$ is important in rank-based statistics where it appears, e.g., 
in many efficacy expressions (Gibbons and Chakraborti, 1992) and in the asymptotic variance of the Hodges-Lehmann estimator.

Källberg et al.\ (2013) study $U$-statistic estimators of $q_{2,0}$ and $h_2$ based on the number of $\ep$-close observations in a sample from a stationary $m$-dependent sequence, 
and establish conditions for their consistency and weak convergence.
In this paper, we investigate the corresponding $U$-statistic estimators for a more general class of functionals, including the (two-density) functional $q_{1,1}$, 
and give alternative conditions for consistency.
Additionally, we propose estimators for $q_{k,l}$ based on incomplete $U$-statistics and show that they are suitable 
when no assumptions, except for $m$-dependence, can be made about the dependencies among the observations.
Fan and Ullah (1999) obtain consistency of kernel estimates of $q_{k,l}$ under dependence weaker than $m$-dependence, 
but do not provide rates of convergence. 
Wavelet-methods for estimation of $q_{2,0}$ with dependent observations are considered in Hosseinioun et al.\ (2009) and Chesneau et al.\ (2013). 

There are numerous studies about nonparametric estimation of $q_{2,0}$ and related functionals for independent observations, and typical methods include those based on nearest neighbors 
(Leonenko et al., 2008), kernels (Bickel and Ritov, 1988, Gine and Nickl, 2008, Chac\'{o}n and Tenreiro, 2012), and orthogonal projection (Laurent, 1997, Tchetgen et al., 2008). Kernel estimators of $q_{2,0}$ (and other integral functionals) from a contaminated sample are studied in Delaigle and Gijbels (2002).

Note that we can introduce  'plug-in' estimators $\hat{D}_n$ and $\hat{h}_n$ of the functionals $D=q_{2,0} -q_{1,1} +q_{2,0}$ and $h_2=-\log(q_{0,2})$, respectively, by replacing $\{q_{k,l}\}$ with their estimators in the corresponding expressions.
Then it straightforward to combine the results in this paper with the conventional limit theory for plug-in estimators to derive consistency properties of the estimators $\hat{D}_n$ and $\hat{h}_n$. 
For applications of such nonparametric divergence and entropy estimators, we refer to, e.g., Leonenko et al.\ (2008), Leonenko and Seleznjev (2010), Källberg et al.\ (2013), and references therein.

First we introduce some notation.
A sequence $\{W_i\}_{i=1}^\infty$ of random $d$-vectors is said to be $m$-dependent if $\{W_{1},\ldots ,W_{b} \}$ and $\{W_{b+k},\ldots,\ldots,W_{b+k+l}\}, b,k,l =1,2,\ldots$, are independent sets  when $k> m$.
Let $L_a(R^d),a \geq 1$, denote the usual Lebesgue space of real-valued functions in $R^d$.
Define $X$ and $Y$ to be independent random vectors with densities $p_X(x)$ and $p_Y(x), x \in R^d$, respectively. 
For $x,y \in R^d$, let $d(x,y):=|x-y|$ denote the Euclidean distance in $R^{d}$. 
Throughout the paper, we assume that $0<\ep =\ep(n) \rightarrow0\mbox{ as }n\rightarrow \infty$.
Let $B_{\ep}(x):=\{y \in R^d:d(x,y)\leq \ep \}, x \in R^d$, be an $\ep$-ball in $R^{d}$ 
with center at $x$ and radius $\ep$. Denote the volume of $B_\ep(x)$ by $\bep =b_1(d)\ep^d$ and 
define the $\ep$-ball probability as $p_{X,\ep}(x) := P(X \in B_\ep(x))$.
Let $I(A)$ be the indicator function of an event $A$. Mean square convergence is denoted by $\ms$.

To evaluate the rate of mean square convergence of our considered estimators, 
we need a smoothness condition for the marginal densities $p_X$ and $p_Y$ (cf.\ Leonenko and Seleznjev, 2010).
Let $H_2^{(\alpha)}(K), 0<\alpha \leq 1, K>0$, be a linear space 
of functions in $L_2(R^d)$ that satisfy an $\alpha$-Hölder smoothness condition in $L_2$-norm with constant $K$, 
i.e., for $p \in H_2^{(\alpha)}(K)$,
$$
\left(\int_{R^d} (p(x-h)-p(x))^2dx \right )^{1/2} \leq K |h|^\alpha, \quad  h \in B_1(0).
$$


The remaining part of this paper is organized as follows. 
In Section 2, we introduce the $U$-statistic estimators of $q_{k,l}$ and present their mean square consistency properties. Some numerical experiments illustrating the obtained rates of convergence are also provided.
In Section 3, we make some conclusions and discuss possible extensions of our approach.
Section 4 contains the proofs of the results in Section 2.
\section{Results}
Without loss of generality, we only present the estimators and the corresponding results for the functionals $q_{2,0}$ and $q_{1,1}$. 

We introduce the following (dependence) conditions for the stationary $m$-dependent sequences $\{X_i\}$ and $\{Y_i\}$. 
Note that $\mc A_1$ and $\mc A_2$ allow dependence between $\{X_i\}$ and $\{Y_i\}$.

\begin{enumerate}
\item[$\mc A_1$] 
The $(2d)$-sequence $\{(X_i,Y_i) \}_{i=1}^\infty$ is strictly stationary and $m$-dependent. 
\item[$\mc A_2$]  For $s,t=1,2,\ldots,$ the random $d$-vectors 
 $X_{1}-X_{1+s}$, and $X_{s}-Y_{t}$ have bounded densities 
 $f^X_{s}(x)$ and $f^{X,Y}_{s,t}(x)$, respectively. 
\end{enumerate}
\textbf{Example 1.}\ (i) $\mc A_1$ and $\mc A_2$ hold if $\{(X_i,Y_i)\}$ is a Gaussian stationary $m$-dependent sequence.\smallskip

\noindent (ii) Under independence, i.e., $\{X_i\}$ and $\{Y_i\}$ are independent and $m=0$, the (minimal) condition $p_X, p_Y \in L_2(R^d)$ is sufficient for $\mc A_2$. \smallskip

\subsection{Estimators based on $U$-statistics}
In this section, we consider the $U$-statistic estimators studied for independent observations in Källberg and Seleznjev (2012). The estimator of $q_{k,l}$, based on equally sized samples $\{X_1,\ldots,X_n\}$ and $\{Y_1,\ldots,Y_n\}$, is defined as $\Qt_{k,l,n} := Q_{k,l,n}/\bep$, where
\begin{align*} 
Q_{2,0,n}  := \frac{1}{\binom{n}{2}} \sum_{1 \leq i<j\leq n}I(d(X_i,X_j)\leq\ep), \quad
Q_{1,1,n} := \frac{1}{n^2} \sum _{i=1}^{n}\sum _{j=1}^{n} I(d(X_i,Y_j) \leq\ep).
\end{align*}
We briefly explain the idea behind this estimator.
For example, in view of the representation $q_{1,1} = \Ex( p_X(Y))$ and the well-known approximation $\bep^{-1} p_{X,\ep}(x)$ of $p_X(x)$, 
an approximation of $q_{1,1}$ can be written as $\qt_{1,1,\ep} := \bep^{-1} \Ex (p_{X,\ep}(Y)) =
 \bep^{-1} P(d(X,Y)\leq\ep)$. 
 Now, the $U$-statistic $Q_{1,1,n}$ may, despite the $m$-dependence, serve as an estimate of the 
 $\ep$-{\it coincidence probability} $P(d(X,Y)\leq\ep)$ for independent observations. 
Hence, $\Qt_{1,1,n}$ is suitable as an estimator of $\qt_{1,1,\ep}$ and thus also of $q_{1,1}$. 
The same reasoning holds for $\Qt_{2,0,n}$.
 \medskip

\noindent \textbf{Remark 1.} 
The assumption of equal sample sizes in this paper is merely for technical convenience, and we claim that the developed method can be used when the sample sizes, say $n_1$ and $n_2$, are unequal and satisfy $n_1/n_2 \to \rho$, for some $\rho > 0$.\medskip

First we study the mean square consistency of $\Qt_{k,l,n}$ under dependence conditions $\mc A_1$ and $\mc A_2$. For the exact asymptotics, we introduce some variance characteristics:
let $g(X_i,Y_i) := \frac{1}{2}(p_Y(X_i) + p_X(Y_i)), i =1,\ldots,n$, and define
\begin{align}\label{sigma2kl}
\sigma_{2,0}^2 & := \Var(p_X(X_1))  + 2\sum_{h=1}^m \Cov(p_X(X_1), p_X(X_{1+h}), \\
\sigma_{1,1}^2 & := \Var(g(X_1,Y_1))+2\sum_{h=1}^m \Cov(g(X_1,Y_1),g(X_{1+h},Y_{1+h})). \notag
\end{align}

The convergence rates in the following theorem are obtained for $\Qt_{k,l,n}$ based on independent observations in Källberg and Seleznjev (2012). Let $L(n), n \geq 1$, be a slowly varying function satisfying $L(n) \to \infty$ as $n \to \infty$ (e.g., $L(n) = \log(n)$).
\begin{theorem}\label{th:main1}
Let $p_X,p_Y \in L_3(R^d)$ and $\mc A_1$ and $\mc A_2$ hold.
\begin{itemize}
\item[(i)] If $n^2\ep^d \to \infty$, then
\begin{equation*}
\Qt_{k,l,n}  \ms  q_{k,l}  \ninf.
\end{equation*}
\item[(ii)] If $p_X,p_Y \in H_2^{(\alpha)}(K),0<\alpha \leq d/4,$ and $\ep \sim cn^{-2/(4\alpha + d)}, c>0$, then
\begin{equation*}
\Ex(\Qt_{k,l,n}  - q_{k,l})^2 = \Or (n^{-8\alpha/(4\alpha+d)}) \ninf.
\end{equation*}
\item[(iii)] If $p_X,p_Y \in H_2^{(\alpha)}(K),\alpha > d/4$, and $\ep \sim L(n)n^{-1/d}$, then
\begin{equation*}
\Ex(\Qt_{k,l,n} - q_{k,l})^2 = 4 \sigma^2_{k,l}n^{-1} + \oor(n^{-1}) \ninf.
\end{equation*}
\end{itemize}
\end{theorem}\smallskip
\noindent  \textbf{Remark 2.}\ (i)
The rates in Theorem \ref{th:main1} are optimal with respect to the upper bound for the mean square error in
\begin{equation*}\label{upper}
\Ex(\Qt_{k,l,n}-q_{k,l})^2 \leq \frac{C_1}{n} + \frac{C_2}{n^2\ep^d} + C_3 \ep^{4\alpha}, \quad C_1, C_2, C_3 >0,
\end{equation*}
which follows from the proof of Theorem \ref{th:main1} (see Källberg and Seleznjev, 2012, for a discussion). 
\smallskip

\noindent (ii)\ Källberg et al.\ (2013) show that the one-sample estimator $\Qt_{2,0,n}$ attains the rates in Theorem \ref{th:main1} under a dependence condition for $\{X_i\}$ which is nonequivalent to $\mc A_2$.\medskip

\noindent (iii) In the one-dimensional independent case, Bickel and Ritov (1988) 
study a kernel estimator of $q_{2,0}$ for a H\"{o}lder class of densities contained in $H_2^{(\alpha)}(K)$, and obtain the same rates of convergence as in Theorem \ref{th:main1} (cf.\ also Gin\'{e} and Nickl, 2008, Tchetgen et al., 2008).\medskip

Next, we consider the more general setting where no assumptions, except for $m$-dependence, can be made about the dependencies for $\{X_i\}$ and $\{Y_i\}$, i.e., conditions like $\mc A_2$ are not allowed. 
There is a large class of stationary $m$-dependent sequences not satisfying $\mc A_2$. 
In Example 2, this class is illustrated by three $1$-dependent representatives. Other examples can be obtained by similar constructions.\medskip

\noindent
\textbf{Example 2.} 
\noindent (i) Let $\{Z_i\}_{i=-\infty}^\infty$ be independent $N(0,1)$-variables and consider the sequence
\begin{equation*}
X_t = Z_{t} Z_{t-1}, \quad t \geq 1.
\end{equation*}
The density of $X_{1} - X_{2}$ is given by $f_1^{X}(x) = K_0(|x|/\sqrt{2})/(\pi\sqrt{2})$, where $K_0(\cdot)$ is the modified Bessel function of the second kind (Craig, 1936). The function $K_0(\cdot)$ is unbounded at the origin, implying that $\mc A_2$ is not valid in this case.\medskip

\noindent  (ii) Let $\{U_i\}_{i=-\infty}^\infty$ be a sequence of independent identically distributed (i.i.d.)\ random variables with density $p(x)$ and distribution function $F(x), x \in R^1$, and let
\begin{equation*}
X_t := \max(U_t,U_{t-1}), \quad t \geq 1.
\end{equation*}
The marginal density of $\{X_i\}$ is $p_X(x) = 2F(x)p(x)$. Assumption $\mc A_2$ does not hold since $P(X_1 = X_2) =1/3$.\medskip

\noindent (iii) Let $\{X^*_i\}_{i=1}^\infty$ be an i.i.d.\ sequence with density $p_X(x), x \in R^d$, and define $\{\xi_i\}_{i=1}^\infty$ to be sequence of independent symmetric Bernoulli variables. Assume further that $\{X^*_i\}$ and $\{\xi_i\}$ are mutually independent, and define 
\begin{equation*}\label{ber}
X_t := X^*_{t + \xi_t}, \quad t \geq 1.
\end{equation*}
Then $\{X_i\}$ has marginal density $p_X(x)$. Here $P(X_1 = X_2) = 1/4$ and hence $\mc A_2$ is not satisfied.\medskip 

We prove the following result for $\Qt_{k,l,n}$ without involving dependence conditions like $\mc A_2$.
Note that, for density smoothness $\alpha < d/2$, the obtained rates of convergence are slower than in Theorem 1. 
\begin{theorem}\label{th:main2} Let $p_X, p_Y \in L_3(R^d)$ and $\mc A_1$ hold. 
\begin{itemize}
\item[(i)] If $n\ep^d \to \infty$, then 
$$
\Qt_{k,l,n} \ms q_{k,l} \ninf. 
$$
\item[(ii)] If $p_X,p_Y \in H_2^{(\alpha)}(K), 0< \alpha \leq d/2$, and $\ep \sim cn^{-1/(2\alpha+d)}, c > 0$, then
\begin{equation*}
\Ex(\Qt_{k,l,n} - q_{k,l})^2 = \Or(n^{-4\alpha/(2\alpha+d)}) \ninf.
\end{equation*}
\item[(iii)] If $d=1$, $p_X,p_Y \in H_2^{(\alpha)}(K), \alpha > 1/2$, and $\ep \sim L(n)n^{-1/2}$, then
\end{itemize}
\begin{equation*}
\Ex(\Qt_{k,l,n} - q_{k,l})^2 = 4 \sigma^2_{k,l}n^{-1} + \oor(n^{-1}) \ninf.
\end{equation*}
\end{theorem}\medskip

\noindent {\bf Remark 3.} From the proof of Theorem \ref{th:main2}, we have
\begin{equation}\label{upper2}
\Ex(\Qt_{k,l,n} - q_{k,l})^2 \leq \frac{C_1}{n} + \frac{C_2}{n^2\ep^{2d}} + C_3 \ep^{4\alpha}, \quad 
\end{equation}
for some $C_1,C_2,C_3 > 0$. As in Theorem \ref{th:main1} (cf.\ Remark 2(i)), the rates in Theorem \ref{th:main2} are optimal with respect to the obtained upper bound for the mean square error in \eqref{upper2}. 

\subsection{Estimators based on incomplete $U$-statistics}
In view of Theorem 2 and the corresponding result for $\Qt_{2,0,n}$ in Källberg et al.\ (2013) (see Remark 2(ii)), 
it is natural to ask whether a dependence condition like $\mc A_2$ is necessary for an estimator of $q_{k,l}$ to attain the convergence rates in Theorem 1. Next we show that this is not the case by introducing a modified version of $\Qt_{k,l,n}$.

Let $\{ m_n \}_{n \geq 1} $ be a non-random sequence of positive integers such that $m_n \geq m$ for all $n$ large enough and $m_n/n \to 0 \ninf$. 
For example, let $m_n = \oor(n) \to \infty$ or if an upper bound $M$ of $m$ is known, let $m_n = M$. 
Now consider the index sets
\begin{align*}
 \mc I_{2,0,n} &:= \{(i,j): 1 \leq i < j \leq n, \; j-i > m_n \}, \\
\mc I_{1,1,n} &:= \{(i,j): 1 \leq i , j \leq n, \; |j-i| > m_n \},
\end{align*}
which have cardinalities $\binom{n-m_n}{2}$ and $2\binom{n-m_n}{2}$, respectively. 
We propose the estimator $\Qt^*_{k,l,n} := Q^*_{k,l,n}/\bep$ of $q_{k,l}$, where
\begin{align*}
Q^*_{2,0,n}=Q^*_{2,0,n}(m_n)  := \frac{1}{\binom{n-m_n}{2}}\sum_{(i,j) \in \mc I_{2,0,n}} \!\!\!\! I(d(X_i,X_j)\leq\ep), \\
Q^*_{1,1,n}=Q^*_{1,1,n}(m_n)  := \frac{1}{2\binom{n-m_n}{2}}\sum_{(i,j) \in \mc I_{1,1,n}} \!\!\!\! I(d(X_i,Y_j)\leq\ep).
\end{align*}
These types of reduced forms are usually referred to as {\it incomplete U-statistics} (see, e.g., Lee, 1990). 
Here, the idea is to alleviate the impact by $m$-dependence by only including terms with $|i-j|$ relatively large. Note that for $n$ large enough the expectation of $\Qt^*_{k,l,n}$ is the same as if the observations were independent.
The question is how the variance behaves.

In Theorem \ref{th:inc}, we show that the incomplete estimator $\Qt^*_{k,l,n}$ satisfies the properties  for $\Qt_{k,l,n}$ in Theorem 1, but without applying $\mc A_2$ or any other similar dependence condition.
As before, $L(n)$ is a slowly varying function and $L(n) \to \infty \ninf$.

\begin{theorem}\label{th:inc}
Let $p_X,p_Y \in L_3(R^d)$ and $\mc A_1$ hold.
\begin{itemize}
\item[(i)] If $n^2\ep^d \to \infty$, then
\begin{equation*}
\Qt^*_{k,l,n}  \ms  q_{k,l}  \ninf.
\end{equation*}
\item[(ii)] If $p_X,p_Y \in H_2^{(\alpha)}(K),0<\alpha \leq d/4,$ and $\ep \sim cn^{-2/(4\alpha + d)}, c>0$, then
\begin{equation*}
\Ex(\Qt^*_{k,l,n}  - q_{k,l})^2 = \Or (n^{-8\alpha/(4\alpha+d)}) \ninf.
\end{equation*}
\item[(iii)] If $p_X,p_Y \in H_2^{(\alpha)}(K),\alpha > d/4$, and $\ep \sim L(n)n^{-1/d}$, then
\begin{equation*}
\Ex(\Qt^*_{k,l,n} - q_{k,l})^2 = 4 \sigma^2_{k,l}n^{-1} + \oor(n^{-1}) \ninf.
\end{equation*}
\end{itemize}
\end{theorem}

\subsection{Numerical experiments} 
In this section, we study our estimation method for some one-dimensional sequences with density smoothness $\alpha=1$. Monte-Carlo estimation of the mean square error (based on 5000 simulations) is used to evaluate the finite sample performance of the proposed estimators, 
and to empirically verify their rate of convergence for some sample sizes between $n=100$ and $n=1000$.\medskip 

\begin{figure}[htbp]
\begin{center}
\includegraphics[width=\textwidth]{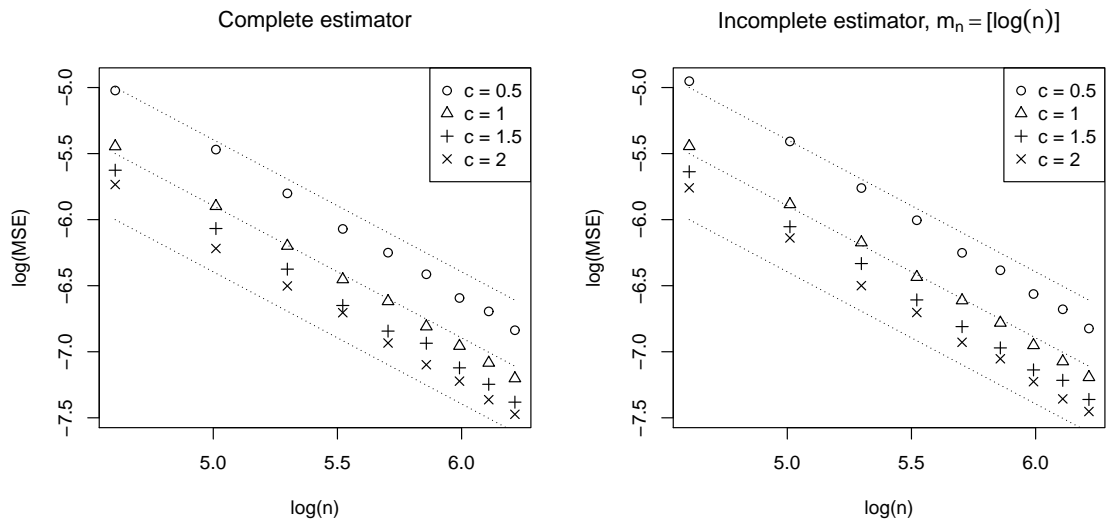}
\caption{Behavior of plug-in estimators, $\hat{D}_n$ and $\hat{D}^*_n$, of the divergence $D$ between the marginal distributions, $N(0,1)$ and $N(1,3/4)$, of  2-dependent Gaussian  sequences. 
Sample sizes $n$ between $100$ and $1000$, and $\ep = c \log(n)n^{-1}$, for certain $c>0$.
Log-log plots of the estimated mean square error (MSE), based on 5000 simulations, against $n$ for $\hat{D}_{n}$ (left) and $\hat{D}^*_n$ (right), 
where the dashed straight lines have negative unit slope.} 
\label{fig1}
\end{center}
\end{figure}

\noindent \textbf{Example 1.} We consider plug-in estimators for the introduced quadratic divergence 
\begin{equation*}
D = \int_{R^d}(p_X(x)-p_Y(x))^2\, dx =q_{2,0} - 2q_{1,1} +  q_{0,2},
\end{equation*}
defined as $\hat{D}_{n} := \Qt_{2,0,n}-2\Qt_{1,1,n}+\Qt_{0,2,n}$ and 
$\hat{D}^*_{n} := \Qt^*_{2,0,n}-2\Qt^*_{1,1,n}+\Qt^*_{0,2,n}$, 
where $\Qt_{2,0,n}$ corresponds to $\Qt_{0,2,n}$ based on the sample $\{Y_1,\ldots, Y_n\}$ and similar for $\Qt^*_{0,2,n}$.
With $[\cdot]$ denoting integer part, we let $m_n:= [\log(n)]$ for the incomplete estimators $\{\Qt^*_{k,l,n}\}$.  
Samples are taken from the $2$-dependent Gaussian moving average time series
\begin{equation*}
X_t  := \frac{1}{\sqrt{3}}(Z_{t} + Z_{t-1}+Z_{t-2}) \quad \mbox{and} \quad Y_t  := \frac{1}{2}(W_{t} - W_{t-1} +  W_{t-2}) +1, \quad n \geq 1,
\end{equation*}
where $\{Z_i\}_{i=-\infty}^\infty$ and $\{W_i\}_{i=-\infty}^\infty$ are independent sequences of i.i.d.\ standard normal variables. 
We have marginal distributions $X_1 \sim N(0,1)$ and $Y_1 \sim N(1,3/4)$, which implies $D\approx 0.155$. The performance of the estimators $\hat{D}_n$ and $\hat{D}^*_n$ is evaluated with $\ep = c\log(n)n^{-1}$ for different values of 
$c>0$. 
Note that Theorems 1 and 3 imply that $\hat{D}_{n}$ and $\hat{D}^*_{n}$ attain the $n$ rate of mean square convergence for these asymptotics of $\ep$. The log-log plots of estimated mean square error against sample size and the straight lines (with negative unit slope) in Figure~\ref{fig1} illustrate this rate of convergence for these values of $n$, $m_n$, and $\ep$. Moreover, the difference in mean square performance between the estimators seems to be negligible. 
\bigskip

\begin{figure}[htbp]
\begin{center}
\includegraphics[width=\textwidth]{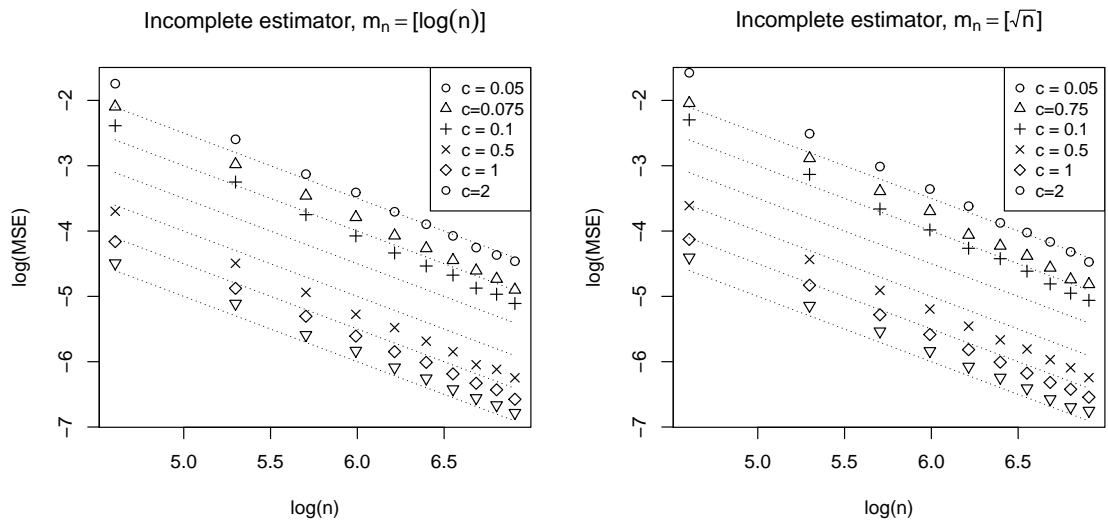}
\caption{{
Behavior of the incomplete estimator $\Qt^*_{2,0,n}$ when $\{X_i\}$ is $2$-dependent with $p_X(x) = e^{-x}, x>0$. 
Sample sizes $n$ between $100$ and $1000$, and $\ep = c \log(n)n^{-1}$, for some $c>0$. 
Log-log plots of estimated mean square error (MSE) based on 5000 simulations 
against $n$, where the dashed straight lines have negative unit slope.
Two versions of $\Qt^*_{2,0,n}$, with $m_n = [\log(n)]$ (left) and $m_n=[\sqrt{n}]$ (right).}}
\label{fig2}
\end{center}
\end{figure}

\noindent \textbf{Example 2.} 
We consider estimation of $q_{2,0}$ for the stationary $2$-dependent sequence 
\[
X_t := \min(Z_{t-2},Z_{t-1},Z_{t}), \quad t \geq 1,
\]
where $\{Z_i\}_{i=-\infty}^\infty$ are independent and exponentially distributed, that is, $Z_i \sim Exp(1/3)$. In this case, the marginal distribution is exponential, $X_t \sim Exp(1), t \geq 1$, and hence $q_{2,0} = 1/2$. We consider two variants of the incomplete estimator $\Qt^*_{2,0,n}$, with sequences $m_n = [\log(n)]$ and $m_n=[\sqrt{n}]$, respectively, and $\ep = c\log(n)n^{-1}$ for some values of $c>0$. Theorem \ref{th:inc} gives that the estimators converge at the rate $n$ under these conditions.  
The log-log plots in Figure \ref{fig2} indicate that this is reasonable for both estimators. Furthermore, no essential differences in mean square performance between the estimators are revealed by these simulations. \medskip 


\section{Concluding remarks and further research}
In this paper, we introduce conditions for the stationary $m$-dependent sequences $\{X_i\}$ and $\{Y_i\}$ under which the estimator $\Qt_{k,l,n}$ attains the same rate of mean square convergence as obtained under independence in Källberg and Seleznjev (2012). Additionally, an incomplete version of $\Qt_{k,l,n}$ is proposed and shown to converge at these rates even without imposing dependence assumptions (beyond $m$-dependence) on $\{X_i\}$ and $\{Y_i\}$. Possible future research is to study these estimators for a wider class of stationary sequences.

It would be of interest to investigate the weak convergence (e.g., asymptotic normality and Poisson convergence) of the considered estimators. The incomplete estimator $\Qt^*_{k,l,n}$ is especially promising; its wide applicability with respect to consistency should, at least to some extent, hold for other asymptotic properties. In fact, it seems that Hoeffding's projection method for asymptotic normality (Sen, 1963) can be used, and we make the following conjecture: given the assumptions in Theorem \ref{th:inc}{\it (iii)}, and if $\sigma^2_{k,l} >0$, then
\begin{equation}\label{norm}
\sqrt{n}(\Qt^*_{k,l,n} - q_{k,l}) \weak N(0,4\sigma^2_{k,l}) \ninf,
\end{equation}
where $\weak$ denotes convergence in distribution.
Because no dependence assumptions except for $m$-dependence are required, \eqref{norm} would, in view of the corresponding result for $\Qt_{2,0,n}$ in Källberg et al.\ (2013), be a significant contribution. 

An important practical problem is how to find a suitable 'bandwidth' $\ep$ for a given sample size. 
Data-driven choice of $\ep$ is an important research question to be addressed in future studies.
Moreover, for the incomplete estimator, we have to make a proper choice of $m_n$, which can be difficult if there is little or no information about $m$. 

\section{Proofs}
\begin{lemma}[Ch.\ 2, Lee, 1990]\label{lemma:lee} 
There are $\binom{n-m}{2}$ pairs of integers $1 \leq i < j \leq n$ that satisfy $j-i >m$.
\end{lemma}

\begin{lemma}\label{lemma2} 
For $j=1,\ldots, 4$, 
assume that the random vector $U_j$ has density $p_{U_j} \in L_3(R^d)$, 
and let $\ph_{U_j}(x) := \bep^{-1} p_{U_j,\ep}(x)$. Then
\[
\Ex \ph_{U_1,\ep}(U_2) \to \Ex p_{U_1}(U_2)  \quad \mbox{and} \quad 
\Ex( \ph_{U_1,\ep}(U_2)\ph_{U_3,\ep}(U_4)) \to \Ex( p_{U_1}(U_2) p_{U_3}(U_4) )\mbox{ as } \ep \to 0,
\]
and hence $\Cov (\ph_{U_1,\ep}(U_2),\ph_{U_3,\ep}(U_4)) \to \Cov (p_{U_1}(U_2),p_{U_3}(U_4))$ as $\ep \to 0$.
\end{lemma}
\medskip

\noindent
The proof of Lemma \ref{lemma2} follows the same steps as that of Lemma 3 in Källberg et al.\ (2013) and therefore  is omitted.
\begin{lemma}\label{lemma3}
If $\mc A_1$ and  $\mc A_2$ hold, then \, $\max_{s, t} P(d(X_s,Y_t)\leq\ep)$\, and\, $\max_{s\neq t} P(d(X_s,X_t)\leq\ep)$ are both $\Or(\ep^d) \mbox{ as } \ep \to 0$.
\end{lemma}
\medskip
\noindent {\it Proof:} 
Let $f^{X,Y}_{s,t}(x)$ be the density of $X_s-Y_{t}$. 
From the boundedness of $f^{X,Y}_{s,t}(x)$, we get
\begin{align}\label{P(Xi,Yj)}
\bep^{-1}P(d(X_s,Y_{t})\leq\ep) &= \bep^{-1}P(|X_s-Y_{t}|\leq\ep) =\bep^{-1}\int_{|x|\leq\ep} f^{X,Y}_{s,t}(x) dx \\  
&  \leq \sup_{x} f^{X,Y}_{s,t}(x) < \infty. \notag 
\end{align}
Moreover, the $m$-dependence and stationarity of $\{(X_i,Y_i) \}$ give that 
$P(d(X_s,Y_t)\leq\ep)$ take on a finite number of values as $\{s,t\}$ vary and so the assertion about $\max P(d(X_s,Y_t)\leq\ep)$ is implied by \eqref{P(Xi,Yj)}. A similar argument applies for the other case. The proof is complete.\qed \\\\

\noindent In the following, we only provide proofs for the two-sample case $(k,l) = (1,1)$, since those for
$\Qt_{2,0,n}$ and $\Qt^*_{2,0,n}$ are similar. As before, $X$ and $Y$ are independent random vectors with densities $p_X(x)$ and $p_Y(x),x \in R^d$, respectively. Note that Lemma \ref{lemma2} gives
\begin{equation}\label{qt11}
\qt_{1,1,\ep}:=\bep^{-1}P(d(X,Y)\leq\ep) =\Ex p_{X,\ep}(Y) \to q_{1,1} \mbox{ as } \ep \to 0.
\end{equation}\medskip

\noindent
{\it Proof of Theorem \ref{th:main1}. }{\it (i)} The assertion is proved by showing that $\Ex \Qt_{1,1,n} \to q_{1,1}$ and $\Var(\Qt_{1,1,n}) \to 0$.
The argument is similar to that of the corresponding result for $\Qt_{2,0,n}$ in K\"{a}llberg et al.\ (2013). 

First we consider $\Ex \Qt_{1,1,n}$.
By applying Lemma \ref{lemma:lee} and the $m$-dependence of $\{(X_i,Y_i)\}$, it follows that there are $2\binom{n-m}{2}$ summands in $\sum_{i,j} I(d(X_i,Y_j)\leq\ep)$ where $X_i$ and $Y_j$ are independent. Hence, 
\begin{align}\label{sum}
\Ex \Qt_{1,1,n}  & = \frac{1}{n^2 \bep} \sum_{1\leq i,j \leq n} P(d(X_i,Y_j) < \ep) \\
& =  \frac{2\binom{n-m}{2}}{n^2 \bep}  P(d(X,Y)\leq\ep) +
\frac{1}{n^2 \bep} \sum_{\substack{1 \leq i,j \leq n\\ |i-j| \leq m }} P(d(X_i,Y_j)\leq\ep). \notag
\end{align}
The last sum in \eqref{sum} consists of $n^2 - 2 \binom{n-m}{2} = \Or(n)$ terms. 
Consequently, from Lemma \ref{lemma3} and \eqref{qt11}, we get
\begin{equation}\label{exp}
\Ex \Qt_{1,1,n} = \qt_{1,1,\ep} + \Or \left( \frac{1}{n} \right) \to q_{1,1} \ninf.
\end{equation}

Next we study the variance of $\Qt_{1,1,n}$.
For $Z_i :=(X_i,Y_i), i =1, \ldots,n$, let
$$
f_\ep(Z_i, Z_j)  := I(d(X_i,Y_j)\leq\ep) + I(d(X_j,Y_i)\leq\ep),
$$
and consider the following reduced form of $N_n:= \sum_{i,j} I(d(X_i,Y_j)\leq\ep)$:
\begin{equation}\label{N*}
N^*_n  := \sum_{\substack{1 \leq i,j \leq n \\ |j-i| > m}} I(d(X_i,Y_j)\leq\ep) = \sum _{\substack{1 \leq i<j \leq n \\ j-i > m}} f_\ep(Z_i,Z_j). 
\end{equation}
%
%
First we study the variance of $N^*_n$:
\begin{equation}\label{cov}
\Var(N_n^*)  =  \sum_{\substack{s_1 < s_2\\  s_2-s_1 >m}}\sum_{\substack{t_1 < t_2\\  t_2-t_1 >m}} \Cov(f_\ep(Z_{s_1},Z_{s_2}), f_\ep(Z_{t_1},Z_{t_2})).
\end{equation}
We need to count the number of terms of different types in this sum.
Let
\begin{equation}\label{f1}
f_\ep^{(1)}(z) : = \Ex f_\ep(z, Z_1)=p_{Y,\ep}(x) + p_{X,\ep}(y), \quad z:=(x,y) \in R^{2d},
\end{equation}
and recall the definition $g(X_i,Y_i) := \frac{1}{2}(p_Y(X_i) + p_Y(Y_i)), i =1,\ldots, n$. 
Now, by using results from Ch.\ 2.4.1.\ in Lee (1990) and Lemma 2, we obtain the following results about the summands in \eqref{cov}:
\begin{itemize}
\item[1)] $|s_i - t_j|>m$ for $i,j = 1,2$. Then covariance \eqref{cov} is zero since all random vectors involved are independent.
\item[2)] One of the differences $|s_i - t_j|$ is zero and the remaining three are greater than $m$. 
For such cases, by conditioning, \eqref{f1}, and Lemma \ref{lemma2}, we get
\begin{align*}
C_{0,\ep}&:=\Cov(f_\ep(Z_{s_1}, Z_{s_2}), f_\ep (Z_{s_1},Z_{t_2})) \\
& =  \Var(f_\ep^{(1)}(Z_{s_1})) \sim 4\bep^{2} \Var(g(X_{s_1},Y_{s_1}))  \ninf.
\end{align*}
The number of these terms is $6\binom{n-2m}{3}$.
\item[3)] For $h=1,\ldots,m$, $0<|s_i-t_j|=h \leq m$ for exactly one of the four differences and the rest are greater than $m$. 
In this case, we obtain from \eqref{f1} and Lemma \ref{lemma2} that
\begin{align*}
C_{h,\ep} &:=\Cov(f_\ep(Z_{s_1}, Z_{s_2}),f_\ep(Z_{s_1+h},Z_{t_2})) = 
\Cov(f_\ep^{(1)}(Z_{s_1}),f_\ep^{(1)}(Z_{s_1+h})) \\ &\sim 4\bep^{2}\Cov(g(X_{s_1},Y_{s_1}),g(X_{s_1+h},Y_{s_1+h}))\ninf.
\end{align*}
There are $12\binom{n-2m-h}{3}$ such terms.
\item[4)] 
From \eqref{qt11}, we get that the remaining $\Or(n^2)$ terms satisfy
\begin{align*}
|\Cov(f_\ep(Z_{s_1},Z_{s_2}),f_\ep(Z_{t_1},Z_{t_2}))| & \leq \Var (f_\ep(Z_{s_1},Z_{s_2})) \leq 2 \Ex f_\ep(Z_{s_1},Z_{s_2}) \\
& = 4P(d(X,Y)\leq\ep)= \Or(\ep^d) \ninf \notag,
\end{align*}
where we have used $s_2-s_1, t_2-t_1 >m$. Consequently, the sum of the terms that are not of type 1)-3) has size $\Or(n^2\ep^{d})$.
\end{itemize}
From these results and the definition of $\sigma^2_{1,1}$ (cf.\ \eqref{sigma2kl}), we obtain
\begin{align}\label{Nn*}
\Var(N_n^*) & =  \sum_{\substack{s_1 < s_2\\  s_2-s_1 >m}}\sum_{\substack{t_1 < t_2\\  t_2-t_1 >m}} \Cov(f_\ep(Z_{s_1},Z_{s_2}), f_\ep(Z_{t_1},Z_{t_2}))\\
& =   6\binom{n-2m}{3}C_{0,\ep} + \sum_{h=1}^m 12\binom{n-2m-h}{3}C_{h,\ep} + \Or(n^2\ep^d) \notag \\
& = 4n^{3}\bep^2\sigma^2_{1,1} + \oor(n^3\ep^{2d}) + \Or( n^2\ep^d ) \ninf. \notag
\end{align}
Moreover, for each 4-tuple $\{s_1,s_2,t_2,t_2 \}$, we get from Lemma \ref{lemma3} that
\begin{align}\label{cov2}
|\Cov(I(d(X_{s_1},Y_{s_2})\leq\ep),I(d(X_{t_1},Y_{t_2})\leq\ep))| &\leq \max_{s,t} \Var(I(d(X_{s},Y_{t})\leq\ep))\\ &= \Or(\ep^d) \ninf, \notag
\end{align}
and hence
\begin{equation}\label{th1:diff}
\Var(N_n-N^*_n) = \Var (\sum_{\substack{1 \leq i,j \leq n \\ |j-i|\leq m}} I(d(X_i,Y_j)\leq\ep) ) = \Or(n^2 \ep^d) \ninf,
\end{equation}
which follows from \eqref{cov2} since $N_n-N^*_n$ consists of $\Or(n)$ terms. Now \eqref{Nn*} and \eqref{th1:diff} give
\begin{align}\label{VQt1}
\Var(\Qt_{1,1,n}) & =n^{-4} \bep^{-2} \Var(N_n) \leq n^{-4} \bep^{-2} 2(\Var(N_n^*) + \Var(N_n-N_n^*)) \\
& = \Or \left(\frac{1}{n} \right) + \Or \left(\frac{1}{n^2\ep^d} \right) \ninf. \notag
\end{align}
Moreover, the condition $n^2 \ep^d \to \infty$ and \eqref{VQt1} yield $\Var(\Qt_{1,1,n}) \to 0$, which together with \eqref{exp} proves the assertion. \\

\noindent {\it (ii)} 
As in Källberg and Seleznjev (2012), the smoothness condition $p_X, p_Y \in H_2^{(\alpha)}(K)$ gives
\begin{equation}\label{e2a}
|\qt_{1,1,\ep} - q_{1,1}| \leq \frac{1}{2}K^2 \ep^{2\alpha}, 
\end{equation}
and so \eqref{exp} yields
\begin{equation}\label{biasBound}
|q_{1,1} - \Ex \Qt_{1,1,n}| \leq \frac{1}{2}K^2 \ep^{2\alpha} + \Or \left( \frac{1}{n} \right) \ninf.
\end{equation}
Further, from the conditions $\ep \sim cn^{-2/(4\alpha + d)}, c>0$, and $0< \alpha \leq d/4$, we get
\begin{equation}\label{lowrate}
\ep^{2\alpha} \sim c^{2\alpha}n^{-4\alpha/(4\alpha + d)} \geq c^{2\alpha} n^{-1/2} \quad\mbox{and} \quad n^2\ep^{d} \sim c^dn^{8\alpha/(4\alpha+d)} \leq n.
\end{equation}
Now the stated rate of mean square convergence follows by combining \eqref{VQt1} and \eqref{biasBound} with \eqref{lowrate}.\\

\noindent {\it (iii) }  First we obtain the order of the bias. The assumptions $\ep \sim L(n)n^{-1/d}$ and $\alpha > d/4$ imply $\Or(\ep^{2\alpha}) = \oor(n^{-1/2})$, 
and hence \eqref{biasBound} yields
\begin{equation}\label{bias}
|q_{1,1} - \Ex \Qt_{1,1,n}|
 = \oor \left( \frac{1}{n^{1/2}} \right) \ninf.
\end{equation}

Next we consider $\Var(\Qt_{1,1,n})$.
First note that limits \eqref{Nn*} and \eqref{th1:diff} together with the assumption $n\ep^d \sim L(n) \to \infty$ give
\begin{align}\label{th1:rates}
n^{-4}\bep^{-2}\Var(N^*_n) & = 4\sigma^2_{1,1}n^{-1} + \oor(n^{-1}),\\
 n^{-4}\bep^{-2}\Var(N_n^*-N_n) & = \oor(n^{-1}) \ninf.  \notag
\end{align}
We get
\begin{align}\label{th1:dec}
\Var(\Qt_{1,1,n}) & = n^{-4}\bep^{-2}\Var(N_n)\\
&= n^{-4}\bep^{-2}(\Var(N_n^*) + \Var(N_n-N_n^*) + 2 \Cov(N^*_n,N_n-N_n^*)) \notag \\
				 & = 4\sigma^2_{1,1}n^{-1} + \oor(n^{-1}) \ninf, \notag
\end{align}
where the last equality follows from \eqref{th1:rates} and by applying the Cauchy-Schwarz inequality to the covariance. 
Now the assertion follows from \eqref{bias} and \eqref{th1:dec}, which completes the proof. \qed \\\\

\noindent In the following proofs of Theorems 2 and 3, we use some intermediate results from the proof of Theorem 1, for example, \eqref{sum}, \eqref{Nn*}, \eqref{e2a}, \eqref{lowrate}, and a weakened version of \eqref{th1:rates},
which is allowed since these hold without assumption $\mc A_2$.\\\\

\noindent
{\it Proof of Theorem \ref{th:main2}}. {\it (i) }  First consider the expectation $\Ex \Qt_{1,1,n}$. 
We use the decomposition of $\Ex \Qt_{1,1,n}$ in \eqref{sum} .
Note that the last sum in \eqref{sum} consists of $\Or(n)$ terms, each bounded by~1. 
Hence,
\begin{align}\label{exp2}
\Ex \Qt_{1,1,n}&  = \frac{2\binom{n-m}{2}}{n^2} \bep^{-1} P(d(X,Y)\leq\ep) + \Or \left( \frac{1}{n \ep^d} \right) \\
&= \qt_{1,1,\ep} + \Or \left( \frac{1}{n \ep^d} \right) \ninf. \notag
\end{align}

Next we study $\Var(\Qt_{1,1,n})$.
The variance of $\Or(n)$ zero-one variables is of size $\Or(n^2)$ and thus
\begin{equation}\label{th2:diff}
\Var(N_n-N^*_n) = \Var (\sum_{\substack{1 \leq i,j \leq n \\ |j-i|\leq m}} I(d(X_i,Y_j)\leq\ep) ) = \Or(n^2) \ninf.
\end{equation}
Then, similarly as in \eqref{VQt1}, we get
\begin{equation}\label{varQ2}
\Var(\Qt_{1,1,n}) = \Or\left( \frac{1}{n} \right) + \Or\left( \frac{1}{n^2\ep^{2d}} \right) \ninf.
\end{equation}
Finally, by combining \eqref{exp2} and \eqref{varQ2} with the assumption $n\ep^d \to \infty$, we get $\Ex \Qt_{1,1,n} \to q_{1,1}$ and $\Var( \Qt_{1,1,n}) \to 0$. The statement follows.
\\

\noindent {\it (ii)} The argument is similar to that of Theorem 1{\it (ii)}.
For the bias, we get from \eqref{e2a} and \eqref{exp2} that
\begin{equation}\label{biasBound2}
|q_{1,1}-\Ex \Qt_{1,1,n}| \leq C\ep^{2\alpha} + \Or\left(\frac{1}{n\ep^d} \right) \ninf.
\end{equation}
Moreover, the assumptions $\ep \sim cn^{-1/(2\alpha + d)}, c>0$, and $0< \alpha \leq d/2$ lead to
\begin{equation}\label{lowrate2}
\ep^{2\alpha} \sim c^{2\alpha}n^{-2\alpha/(2\alpha + d)} \geq c^{2\alpha} n^{-1/2} \quad\mbox{and} \quad n^2\ep^{2d} \sim c^{2d}n^{4\alpha/(2\alpha+d)} \leq c^{2d}n.
\end{equation}
The asserted rate of convergence is implied by \eqref{varQ2}, \eqref{biasBound2}, and \eqref{lowrate2}.\\

\noindent {\it (iii)} 
By applying the conditions $\ep \sim L(n)n^{-1/2}$ and $\alpha > 1/2$, 
we get $\Or(\ep^{2\alpha}) = \oor(n^{-1/2})$ and $n^{1/2}\ep^1 \to \infty$. 
Hence, \eqref{biasBound2} yields
\begin{align}\label{th2:exp}
|q_{1,1} - \Ex \Qt_{1,1,n}| = \oor \left( \frac{1}{n^{1/2}}\right) \ninf.
\end{align}

For $\Var(\Qt_{1,1,n})$, we follow the same steps as in Theorem 1. 
If we use \eqref{th1:rates} and decomposition \eqref{th1:dec}, but with the weaker rate \eqref{th2:diff} for $\Var(N_n - N_n^*)$ in place of \eqref{th1:diff}, it follows that
\begin{align}\label{th2:var}
\Var(\Qt_{1,1,n}) &  = n^{-4}\bep^{-2}(\Var(N_n^*) + \Var(N_n-N_n^*) + 2 \Cov(N^*_n,N_n-N_n^*)) \\
& = 4\sigma^2_{1,1}n^{-1} + \oor(n^{-1}) + \Or\left( \frac{1}{n^2\ep^2} \right) + \Or\left( \frac{1}{n^{3/2}\ep} \right) \notag \\
& = 4\sigma^2_{1,1}n^{-1} + \oor(n^{-1})  \ninf, \notag
\end{align}
where the last equality follows since $n \ep^2 \sim L(n)^2 \to \infty$ by assumption. 
The claimed mean square property of $\Qt_{1,1,n}$ is implied by  \eqref{th2:exp} and \eqref{th2:var}. This completes the proof. \qed \\\\

\noindent
{\it Proof of Theorem \ref{th:inc}}. The proof is straightforward since the behavior of the variable $\sum_{|i-j|>m_n} I(d(X_i,Y_j)\leq\ep)$, 
which is the basis of $\Qt^*_{1,1,n}$, is very similar to that of $N_n^*$ (cf. \eqref{N*}) 
and so, we can use the asymptotics of $N^*$ derived in the proof of Theorem 1.  \\

\noindent {\it (i)} Assume without loss of generality that $m \leq m_n < n$. Then we have a simple expression for the expectation:
\begin{equation}\label{E}
\Ex \Qt^*_{1,1,n} = \bep^{-1} P(d(X,Y)\leq\ep) = \qt_{1,1,\ep},
\end{equation}
which combined with \eqref{qt11} yields $\Ex \Qt^*_{1,1,n} \to q_{1,1}$.

For the variance of $\Qt^*_{1,1,n}$, we consider the variance of
\begin{equation*}
M^*_n  := \sum_{(i,j) \in \mc I_{1,1,n}} I(d(X_i,Y_j)\leq\ep). 
\end{equation*}
Note that formula \eqref{Nn*} for $N^*_n$ can be applied also for $M^*_n$, 
because the number of the remaining terms (of type 4) decreases when $m_n \geq m$ is used instead of $m$. 
Hence, since $m_n/n \to 0$, it follows directly from \eqref{Nn*} that
\begin{align*}
\Var(M_n^*) & =   6\binom{n-2m_n}{3}C_{0,\ep} + \sum_{h=1}^{m} 12\binom{n-2m_n-h}{3}C_{h,\ep} + \Or(n^2\ep^d) \notag \\
& = 4n^{3}\bep^2\sigma^2_{1,1} + \oor(n^3\ep^{2d}) + \Or( n^2\ep^d ) \ninf \notag
\end{align*}
and consequently, we have
\begin{align}\label{V}
\Var(\Qt^*_{1,1,n}) & = \binom{n-m_n}{2}^{-2} \frac{1}{4} \bep^{-2}\Var(M_n^*) \\
& = 
4\sigma^2_{1,1}n^{-1} + \oor(n^{-1}) + \Or \left( \frac{1}{n^2\ep^d} \right) \ninf. \notag
\end{align}
The condition $n^2\ep^d \to \infty$ and \eqref{V} give $\Var(\Qt^*_{1,1,n}) \to 0$. The assertion follows.\\

\noindent {\it (ii) } The order of the bias $|\Ex\Qt^*_{1,1,n}-q_{1,1}|$ follows directly by combining \eqref{e2a} and \eqref{lowrate} with \eqref{E}. 
Similarly, we see that the order of $\Var(\Qt^*_{1,1,n})$ follows from \eqref{lowrate} and \eqref{V}. 
The stated rate for the mean square convergence of $\Qt^*_{1,1,n}$ is thus obtained.\\

\noindent {\it (iii) } The condition $\ep \sim L(n)n^{-1/d}$ yields $n\ep^d \to \infty$, and so \eqref{V} implies $\Var(\Qt^*_{1,1,n}) = 4\sigma^2_{1,1}n^{-1} + \oor(n^{-1})$. For the bias, we get from \eqref{e2a} and \eqref{E} that $|E \Qt^*_{1,1,n} - q_{1,1}| =|\qt_{1,1,\ep} - q_{1,1}| = \Or(\ep^{2\alpha}) = \oor(n^{-1/2})$, where the last equality follows since $\alpha > d/4$ and $\ep \sim L(n)n^{-1/d}$. These asymptotics for the bias and variance imply the claimed mean square property of $\Qt^*_{1,1,n}$. This completes the proof.
\qed
 \\\\
{\noindent {\large \textbf{References}}}

\begin{reflist}

 {\small Bickel, P.J. and Ritov, Y.\ (1988).
 Estimating integrated squared density derivatives: sharp best order of convergence estimates,
{Sankhy{\=a} A} 50, 381--393.}

{\small   Chac{\'o}n, J.\ E., Tenreiro, C.\ (2012). 
Exact and asymptotically optimal bandwidths for kernel estimation of density functionals.
{Methodol.\ Comput.\ Appl.\ Prob.}\ 14, 523--548.}

{\small Chesneau, C., Kachour, M., Navarro, F.\ (2013). 
A note on the adaptive estimation of a quadratic functional from dependent observations.
{\.{I}statistic} 6, 10--26.}

 {\small Craig, C.C.\ (1936). On the frequency function of xy.
{Ann.\ Math.\ Statist}.\ 36, 1--15.}

{\small Delaigle, A., Gijbels, I.\ (2002). Estimation of integrated squared density derivatives from a contaminated sample. {J.\ R.\ Stat.\ Soc.\ Ser.\ B Stat.\ Methodol.}\ 64, 869-886.}

 {\small Fan, J., Ullah, A.\ (1999). 
 On goodness-of-fit tests for weakly dependent processes using kernel method.
{J.\ Nonparametr.\ Stat.}\ 11, 337--360.}

{\small Gibbons, J. D., Chakraborti, S.\ (1992). Nonparametric statistical inference Vol.\ 131. 
Marcel Dekker, New York.}

 {\small Gin\'{e}, E., Nickl, R., (2008). A simple adaptive estimator of the integrated square of a density,           
 {Bernoulli} 14, 47--61.}

{\small Hosseinioun, N., Doosti, H., Niroumand, H.\ A.\ (2009). 
Wavelet-based estimators of the integrated squared density derivatives for mixing sequences. 
{Pakistan J.\ Statist.}\ 25, 341-350.}

{\small Källberg, D., Leonenko, N., Seleznjev, O.\ (2013). 
Statistical estimation of quadratic R\'{e}nyi entropy for a stationary $m$-dependent sequence. 
arXiv preprint arXiv:1303.1743.}

{\small Källberg, D., Seleznjev, O.\ (2012). 
Estimation of entropy-type integral functionals. arXiv preprint arXiv:1209.2544.}

{\small Laurent, B.\ (1997). Estimation of integral functionals of a density and its derivatives.
{Bernoulli} 3, 181--211.}

{\small Lee, A.J.\ (1990). $U$-Statistics: Theory and Practice. Marcel Dekker, New York.}

{\small Leonenko,  N., Pronzato, L., Savani, V.\ (2008).
 A class of R\'{e}nyi information estimators for multidimensional densities.
{Ann.\  Statist.}\ 36, 2153--2182. Corrections (2010). {Ann.\ Statist.}\ 38, 3837--3838.}

{\small  Leonenko, N., Seleznjev, O.\ (2010). 
Statistical inference for the $\ep$-entropy and the quadratic R\'{e}nyi entropy.
{J.\  Multivariate Anal.}\ 101, 1981--1994.}

{\small R\'{e}nyi, A.\ (1970). 
Probability Theory. North-Holland, Amsterdam.}

{\small Sen, P. K.\ (1963). 
On the properties of U-statistics when the observations are not independent. 
{Calcutta Statist.\ Assoc.\ Bull.}\ 12, 69-92.}

{\small Tchetgen, E., Li, L., Robins, J., van der Vaart, A. (2008). 
Minimax estimation of the integral of a power of a density. {Statist.\ Probab.\ Lett.}\ 78, 3307-3311.}



%
%
%

 \end{reflist}
\end{document}